\documentclass[10pt]{amsart}
\usepackage{amsmath, amsthm, amssymb, euscript, enumerate}
\usepackage{pxfonts, graphicx, subfigure, epsfig}
\usepackage{mathtools, ulem, appendix, accents, color}
\usepackage{esvect}

\newcommand{\R}{{\mathbb R}}
\newcommand{\T}{{\mathbb T}}
\newcommand{\Z}{{\mathbb Z}}

\newcommand{\cN}{{\mathcal N}}

\newcommand{\cA}{{\mathcal A}}

\newcommand{\Ss}{{\mathbb S}}

\newcommand{\e}{\varepsilon}
\newcommand{\p}{\partial}
\newcommand{\vp}{\varphi}
\newcommand{\osc}{\operatornamewithlimits{osc}}

\newcommand{\spt}{\operatorname{spt}}

\newcommand{\dist}{\operatorname{dist}}
\newcommand{\diam}{\operatorname{diam}}

\newcommand{\ra}{\rightarrow}

\newcommand\norm[1]{\left\Arrowvert {#1}\right\Arrowvert}

\theoremstyle{plain}
\newtheorem{theorem}{Theorem}[section]

\newtheorem{corollary}[theorem]{Corollary}

\newtheorem{lemma}[theorem]{Lemma}
\newtheorem{proposition}[theorem]{Proposition}

\theoremstyle{definition}
\newtheorem{definition}[theorem]{Definition}
\newtheorem{example}[theorem]{Example}

\theoremstyle{remark}

\newtheorem{remark}[theorem]{Remark}

\numberwithin{equation}{subsection}

\title [Homogenization of the boundary value for the Dirichlet problem]{Homogenization of the boundary value for the Dirichlet problem}

\author{Sunghan Kim}
\address{Department of Mathematical Sciences, Seoul National University, Seoul 08826, Korea }
\email{sunghan290@snu.ac.kr}
\author{Ki-Ahm Lee}
\address{Department of Mathematical Sciences, Seoul National University, Seoul 08826, Korea
\& Korea Institute for Advanced Study, Seoul 02455, Korea}
\email{kiahm@snu.ac.kr}
\author[Henrik Shahgholian ]{Henrik Shahgholian}
\address{Department of Mathematics, Royal Institute of Technology,
  100~44  Stockholm, Sweden}
\email{henriksh@kth.se}
\thanks{S. Kim has been supported by National Research Foundation of Korea (NRF) grant funded by the Korean government (NRF-2014-Fostering Core Leaders of the Future Basic Science Program). K.-A. Lee has been supported by the National Research Foundation of Korea (NRF) grant funded by the Korean government (MSIP) (No. NRF-2015R1A4A1041675). K.-A. Lee also holds a joint appointment with the Research Institute of Mathematics of Seoul National University. H. Shahgholian has been supported in part by Swedish Research Council. This project is part of an STINT (Sweden)-NRF (Korea) research cooperation program. \\
\indent We would like to thank A. Minne and K. Keryan for reading and commenting the first version of the paper, in 2012. Special thanks go to Hayk Aleksanyan for careful reading and extremely valuable comments, specially in regards to some technical aspects as well as existing literature.}

\begin{document}
    
\begin{abstract}
In this paper, we give a mathematically rigorous proof of the averaging behavior of oscillatory surface integrals. Based on ergodic theory, we find a sharp geometric condition which we call irrational direction dense condition, abbreviated as IDDC, under which the averaging takes place. It should be stressed that IDDC does not imply any control on the curvature of the given surface. As an application, we prove homogenization for elliptic systems with Dirichlet boundary data, in $C^1$-domains.
\end{abstract}

\maketitle

\tableofcontents

\section{Introduction}\label{section:intro}

This paper is concerned with periodic homogenization of uniformly elliptic systems accompanied with rapidly oscillating boundary data,
\begin{equation}\tag{$L_\e$}\label{eq:Le}
\begin{dcases}
-\nabla \cdot \left( A\left(\frac{x}{\e}\right) \nabla u_\e(x) \right) = 0 & \text{in }\Omega, \\
u_\e(x) = g \left(x,\frac{x}{\e}\right) & \text{on }\Gamma, 
\end{dcases}
\end{equation}
where $\Omega$ is a bounded domain in $\R^n$ with the boundary $\Gamma$. Under ellipticity and periodicity assumptions on $A$, it is by now classical that the operator $-\nabla\cdot(A(\frac{\cdot}{\e})\nabla)$ is homogenized in the interior of $\Omega$;  see  e.g. \cite{BLP}. In this paper, 
we impose a sharp geometric condition on $\Gamma$, under which the boundary layer homogenization takes place for the system \eqref{eq:Le}. 

Our analysis starts with finding a sharp geometric condition on $\Gamma$, under which one has the convergence,
\begin{equation}\label{eq:osc-int}
\lim_{\e\ra 0}\int_\Gamma g \left(x,\frac{x}{\e}\right) d\sigma_x = \int_\Gamma\int_{\T^n} g(x,y)dyd\sigma_x,
\end{equation}
when $g$ is periodic in the second argument. We set our starting point at studying the oscillatory surface integral \eqref{eq:osc-int}, since the convergence in \eqref{eq:osc-int} is closely related to the homogenization of \eqref{eq:Le} through the Poisson integral representation of the solution $u_\e$. Such a relationship becomes more apparent when the coefficient $A$ is given by a constant mapping.

The study of the surface integral \eqref{eq:osc-int} has been a subject in harmonic analysis as an application of stationary phase analysis, and \eqref{eq:osc-int} is well-known where $\Gamma$ is a boundary of a smooth convex domain, and $g$ is smooth and periodic in its oscillating variable. In fact, equipped with a uniform control over the curvature, one may quantify the rate of convergence in \eqref{eq:osc-int}, as shown in \cite{S}; see also \cite{ASS1} and \cite{ASS2} in this regard. However, justification of the convergence in \eqref{eq:osc-int} becomes more difficult when one loses control on the smoothness and the curvature of $\Gamma$ as well as the smoothness of $g$. 

In this paper, we prove that \eqref{eq:osc-int} holds under irrational direction dense condition (IDDC in the sequel), defined in Definition \ref{definition:iddc}, and that IDDC is indeed a sharp geometric condition which does not imply any curvature control on $\Gamma$. Moreover, this condition can be defined under a very mild regularity assumption on $\Gamma$. The difficulty in handling the oscillatory surface integral in \eqref{eq:osc-int} without any curvature control of $\Gamma$ arises from the non-uniformity of the ergodic property in irrational directions. In addition, the problem becomes more subtle, since we need to select a correct mesoscopic scale, which is small enough so that $\Gamma$ is well approximated by its tangent hyperplanes, but not too small such that those hyperplanes foliate large enough portion of the unit torus. 

First we choose a partition of unity, to localize the integral around a boundary point $z$ having irrational normal direction. The size $r$ of the local neighborhood $Q$ has to be chosen in a way that tangent hyperplane $\Pi$ at $z$ is close enough to $\Gamma$ in $Q$, and that we can ``freeze'' the first argument of $g$ at $z$. However, $r$ should not be taken too small, since we need the domain $\e^{-1}(\Pi\cap Q)$ for the second argument of $g$ to foliate large enough portion of $\T^n$, so that the ergodicity of the irrational normal at $z$ may apply; note that $\e^{-1}(\Pi\cap Q)$ is of size $\e^{-1} r$. Moreover, we need to assure that such an ergodic property implies an averaging behavior of the localized surface integral under our consideration, since we have a weighted integral after localization:  
\begin{equation}\label{eq:osc-wint}
\lim_{\e\ra 0} \left[\fint_{\Pi\cap Q} \psi(x) g\left(z,\frac{x}{\e}\right) d\sigma_x - \fint_{\Pi\cap Q} \psi(x)d\sigma_x\int_{\T^n} g(z,y)dy\right] = 0,
\end{equation}
where $\psi$ is an element of the partition of unity supported in $Q$, with $|\nabla \psi|$ of order $r^{-1}$. Thus, to have \eqref{eq:osc-int}, we need to sort out irrational boundary points for which \eqref{eq:osc-wint} can be made uniform over the mesoscopic scale $r$, which is comparably larger than $\e$ but smaller than $1$, and satisfies all the aforementioned properties. 

In spite of the subtleness in the analysis of \eqref{eq:osc-int}, we overcome such difficulties by means of some classical tools from measure theory and ergodic theory. We prove a weighted version of Weyl's lemma to have \eqref{eq:osc-wint} for each irrational boundary point $z$, and make use of Egoroff's theorem to sort out the irrational directions for which \eqref{eq:osc-wint} can be made uniform. 

Our first main result is stated as follows.

\begin{theorem}\label{theorem:osc-int} Let $\Gamma$ be a compact $C^1$-hypersurface in $\R^n$ ($n\geq 2$) satisfying IDDC and $g\in L^\infty(\Gamma\times\T^n)$ be such that $\{g(z,\cdot)\}_{z\in\Gamma}$ is an equicontinuous family on $\T^n$. Then convergence in \eqref{eq:osc-int} holds. 
\end{theorem}

\begin{remark}\label{remark:osc-int-sharp}
The geometric condition, IDDC, is sharp and does not imply any curvature control. The sharpness is shown in Subsection \ref{subsection:example}. 
\end{remark}

\begin{remark}\label{remark:osc-int-general}
The regularity of $\Gamma$ and $g$ is not sharp  and can be relaxed to a weaker condition.
 We leave this part to interested readers. 
\end{remark}

\begin{remark}\label{remark:osc-int-const}
A direct application of Theorem \ref{theorem:osc-int} is boundary layer homogenization of constant elliptic systems; see Proposition \ref{proposition:const}. The sharpness of IDDC is also true in this framework. On the other hand, we capture some interesting phenomena, when $\Gamma$ does not satisfy IDDC; see Example \ref{example:limit} and Example \ref{example:discont}. 
\end{remark}

Next we prove homogenization of the elliptic system \eqref{eq:Le}, as an application of Theorem \ref{theorem:osc-int}, based on recent development of homogenization of the oscillatory Poisson kernel $P_\e$ for \eqref{eq:Le}; namely, we will use an expansion of the form
\begin{equation}\label{eq:Pe-peri}
P_\e(x,z) = \overline{P}(x,z) w\left(z,\frac{z}{\e}\right) + R_\e(x,z)\quad(x\in\Omega,z\in\Gamma),
\end{equation}
where $\overline{P}$ is the homogenized Poisson kernel, $w$ is the matrix corrector we construct so as to be periodic in its second argument and $R_\e$ is the remainder term decaying uniformly on any $E\times\Gamma$ for $E\Subset\Omega$. It was first proved in \cite{AL2} that 
\begin{equation}\label{eq:Pe-ome}
P_\e(x,z) = \overline{P}(x,z) \omega_\e(z) + R_\e(x,z)\quad (x\in\Omega,z\in\Gamma),
\end{equation}
for some corrector $\omega_\e$, provided $\Gamma$ is $C^{1,\mu}$ ($0<\mu\leq 1$). A quantitative estimate on the remainder term $R_\e$ was established in \cite{KLS}, under a higher regularity assumption on $\Gamma$. Only recently, \eqref{eq:Pe-peri} was shown in \cite{AKMP} with an improved estimate on the remainder term, provided $\Gamma$ is the boundary of a uniformly convex smooth domain. Here we extend the latter to general $C^{1,\mu}$-domains at a cost of losing sharp control on $R_\e$, since then the main argument relies only on techniques irrelevant with the geometry of given boundaries, such as sharp regularity estimates \cite{AL1} and some refined estimates for half-space problems \cite{GM1}, \cite{GM2} and \cite{P}. By means of \eqref{eq:Pe-peri}, we derive from the integral representation formula for \eqref{eq:Le} that 
\begin{equation}\label{eq:ue-lim}
\begin{split}
\lim_{\e\ra 0} u_\e(x) & = \lim_{\e\ra 0}\int_\Gamma \overline{P}(x,z) w\left(z,\frac{z}{\e}\right) g\left(z,\frac{z}{\e}\right) d\sigma_z = \int_\Gamma \overline{P}(x,z) \overline{g}(z) d\sigma_z\quad(x\in\Omega),
\end{split}
\end{equation}
where the effective boundary data $\overline{g}$ is given by 
\begin{equation}\label{eq:gov}
\overline{g}(z) := \int_{\T^n} w(z,y) g(z,y) dy\quad(z\in\Gamma),
\end{equation}
with $w$ defined by \eqref{eq:w}. Here in deriving the second equality in \eqref{eq:ue-lim} we use the convergence result for \eqref{eq:osc-int}, and it is exactly where the geometric condition of $\Gamma$ comes into the play. Now \eqref{eq:ue-lim} shows that, under the sharp geometric condition on $\Gamma$ ensuring \eqref{eq:osc-int}, \eqref{eq:Le} is homogenized to the elliptic system,
\begin{equation}\tag{$\overline{L}$}\label{eq:Lov}
\begin{cases}
-\nabla \cdot( \overline{A} \nabla \overline{u}(x)) = 0 & \text{in }\Omega,\\
\overline{u}(x) = \overline{g}(x) & \text{on }\Gamma,
\end{cases}
\end{equation}
where $\overline{A}$ is the effective coefficient corresponding to $A$. 

Recent development regarding homogenization of \eqref{eq:Le} can be found in \cite{GM1}, \cite{GM2}, \cite{P}, \cite{ASS1}, \cite{ASS2}, \cite{ASS3}, \cite{A}, \cite{AKMP}, \cite{SZ} and \cite{Z}, and the references therein. The existing literature is mainly concerned with the rate of convergence, under some strict geometric condition on the boundary and smoothness of both the boundary and the data. It should be noted that some of the existing results also cover domains which are not necessarily uniformly convex; for instance, see \cite{Z}. Let us also note that a slightly modified argument of \cite{ASS1} enables one to treat domains where at each point of the boundary there is at least one non-varnishing principal curvature. However, this is the first paper, to the best of our knowledge, proving boundary homogenization with a sharp geometric condition on the boundary and under a mild regularity assumption on both the boundary and the boundary data. Nevertheless, let us remark that boundary homogenization of non-divergence type elliptic equations has already been proved under the sharp geometric condition, IDDC, in \cite{CKL} and \cite{CK} for Neumann boundary data and in \cite{F} for Dirichlet boundary data. 

Combining our main result with (enhanced versions of)  already existing results, we can prove the following. 

\begin{theorem}\label{theorem:bdry-hom} Let $A\in C^{1+[\frac{n}{2}]}(\T^n;\R^{n^2\times m^2})$ ($n\geq 2$, $m\geq 1$) be a uniformly elliptic mapping in the sense of \eqref{eq:ellip-A}. Also assume that $\Omega\subset\R^n$ is a bounded $C^{1,\mu}$-domain ($0<\mu\leq 1$) whose boundary $\Gamma$ satisfies IDDC, and $g\in L^\infty(\Gamma\times\T^n;\R^m)$ is such that $\{g(z,\cdot)\}_{z\in\Gamma}$ is an equicontinuous family on $\T^n$. Then the sequence $\{u_\e\}_{\e>0}$ of weak solutions of \eqref{eq:Le} converges, locally uniformly in $\Omega$ as $\e\ra 0$, to the weak solution $\overline{u}$ of \eqref{eq:Lov}.
\end{theorem}

\begin{remark}\label{remark:bdry-hom-A} The $C^{1+[\frac{n}{2}]}$-regularity of $A$ is assumed to allow the Sobolev embedding theory, which will be needed for the regularity of boundary layer correctors. 
\end{remark}

\begin{remark}\label{remark:bdry-hom-C1a-bdry} The $C^{1,\mu}$-regularity of $\Gamma$ is assumed to have the desired approximation \eqref{eq:Pe-peri} of the oscillatory Poisson kernel. As is shown in Proposition \ref{proposition:const}, $C^1$-regularity suffices to have homogenization for constant coefficient case, where the Poisson kernel no longer oscillates in $\e$-scales.  
\end{remark}

\begin{remark}\label{remark:bdry-hom-Linf-gov} The effective boundary data $\overline{g}$ can be discontinuous even though $g$ is assumed to be smooth. This phenomenon accounts for influence of the interior oscillation near the boundary layer. In contrast, the regularity of boundary data is preserved through the limit when interior coefficient is constant (c.f. Remark \ref{remark:const-reg-g}).
\end{remark}

This paper is organized as follows. Basic notation and terminology are introduced in Section \ref{section:notation}. In Section \ref{section:osc-int}, we study averaging behavior of surface integrals of oscillatory functions and prove Theorem \ref{theorem:osc-int}. Section \ref{section:bdry-hom} is devoted to the proof of Theorem \ref{theorem:bdry-hom} which involves analysis of boundary layer correctors. 

%%%%%%%%%%%%%%%%%%%%%%%%%%%%%%%%%%%%%%%%%%%%%%%%%%%%%%%%%%%%%%%%%%%%%%%%%%%%%%%%%%%%%%%%%%%%%%%%
%
%     Notation and Terminology
%
%%%%%%%%%%%%%%%%%%%%%%%%%%%%%%%%%%%%%%%%%%%%%%%%%%%%%%%%%%%%%%%%%%%%%%%%%%%%%%%%%%%%%%%%%%%%%%%%

\section{Notation and Terminology}\label{section:notation}

By $\T^n$ we denote the $n$-dimensional unit torus, and simply write $\T^1$ by $\T$. By $\Ss^{n-1}$ we denote the $(n-1)$-dimensional unit sphere. By $Q_r(z)$ (resp., $B_r(z)$) we denote the $n$-dimensional cube (resp., ball) centered at $z$ with side length (resp., radius) $r$, and simply by $Q_r$ (resp., $B_r$) when $z$ is the origin. We write $x= (x',x_n)\in\R^n$ with $x'\in\R^{n-1}$, and similarly, $Q_r'(z) = Q_r(z) \cap \{x_n=0\}$. Given a vector $e$, by $e_i$ we denote the $i$-th component of $e$. By $\{e^1,\cdots, e^n\}$ we denote the standard basis for $\R^n$ that is $e_i^j=1$ if $i=j$ and $e_i^j=0$ if $i\neq j$. By $[\alpha]$ we denote the greatest integer less than or equal to $\alpha\in\R$. 

The rationality of a direction is defined as follows. 

\begin{definition}\label{definition:vector} 
\begin{enumerate}[(i)]
\item $\nu\in\Ss^{n-1}$ is called a rational direction, if $\nu=\alpha k$ for some $\alpha\in\R$ and $k\in\Z^n$.  
\item $\nu\in\Ss^{n-1}$ is called an irrational direction, if $\nu$ is not a rational direction.
\end{enumerate}
\end{definition}

Now let $\Gamma$ be a $C^1$-hypersurface in $\R^n$ with the inward unit normal mapping $\nu:\Gamma\ra\Ss^{n-1}$. By $\Pi(z)$ we denote by the tangent hyperplane to $\Gamma$ at $z$, 
\begin{equation}\label{eq:Pi}
\Pi(z)   := \{x\in\R^n: \nu(z) \cdot(x-z) = 0\},
\end{equation}
and by $H(z)$ the half-space, 
\begin{equation}\label{eq:H}
H(z) := \{x\in\R^n: \nu(z)\cdot (x-z)>0\},
\end{equation}
whose boundary is exactly the tangent hyperplane $\Pi(z)$. Let us define irrational direction dense condition, abbreviated as IDDC in the sequel, as follows. 

\begin{definition}\label{definition:iddc} We say $\Gamma$ satisfies IDDC, if 
\begin{equation}\label{eq:iddc}
\sigma(\{\text{$z\in\Gamma$ : $\nu(z)$ is rational}\}) = 0,
\end{equation}
where $\nu$ is the unit normal mapping on $\Gamma$ and $\sigma$ is the surface measure on $\Gamma$.
\end{definition} 

Let $A=(A_{ij}^{\alpha\beta})$ with $1\leq i,j\leq n$ and $1\leq \alpha,\beta\leq m$ be a $\R^{n^2\times m^2}$-valued mapping on $\T^n$. We say that $A$ is uniformly elliptic, if there is a constant $0<\lambda<1$ such that
\begin{equation}\label{eq:ellip-A}
\lambda |\xi|^2 \leq A_{ij}^{\alpha\beta}(y)\xi_i^\alpha\xi_j^\beta \leq \frac{1}{\lambda}|\xi|^2\quad\text{for all $\xi\in\R^{n\times m}$ and all $y\in\T^n$}.
\end{equation}
Note that the system \eqref{eq:Le} is written in coordinates as 
\begin{equation*}
\begin{dcases}
- \frac{\p}{\p x_i} \left( A_{ij}^{\alpha\beta} \left(\frac{x}{\e}\right) \frac{\p}{\p x_j} u_\e^\beta(x) \right) = 0 & \text{in }\Omega,\\
u_\e^\alpha (x) = g^\alpha \left(x,\frac{x}{\e}\right) & \text{on }\Gamma,
\end{dcases}
\end{equation*}
for $1\leq \alpha\leq m$. 

%%%%%%%%%%%%%%%%%%%%%%%%%%%%%%%%%%%%%%%%%%%%%%%%%%%%%%%%%%%%%%%%%%%%%%%%%%%%%%%%%%%%%%%%%%%%%%%%
%
%     Oscillatory Surface Integrals
%
%%%%%%%%%%%%%%%%%%%%%%%%%%%%%%%%%%%%%%%%%%%%%%%%%%%%%%%%%%%%%%%%%%%%%%%%%%%%%%%%%%%%%%%%%%%%%%%%

\section{Oscillatory Surface Integrals}\label{section:osc-int}

This section is devoted to the proof of Theorem \ref{theorem:osc-int}. Throughout this section, let $\Gamma$ be a compact $C^1$-hypersurface in $\R^n$ with $n\geq 2$. Also let $\nu:\Gamma\ra\Ss^{n-1}$ be the inward unit normal mapping on $\Gamma$. 

\subsection{A Weighted Version of Weyl's Lemma}\label{subsection:tech}

The main tool is the ergodic theory which links the surface integral to the average integral on the torus, owing to the foliation by hyperplanes with irrational directions. Let us begin with some technical lemmas to build up our strategy. 

\begin{lemma}\label{lemma:weyl}
Let $e'\in \R^{n-1}$ be such that $\frac{(e',1)}{|(e',1)|}\in\Ss^{n-1}$ is an irrational direction, $h\in C(\T)$ and $\phi\in C^1(\R^n)$ be such that $0\leq \phi\leq 1$, $|\nabla\phi|\leq \alpha$, $\phi = 1$ on $Q_1$ and $\spt(\phi)\subset Q_{4/3}$. Then for any $\eta>0$, there exists $\rho_0>0$, which depends only on $\eta$, $n$, $e'$, $\alpha$, $\beta$, $N$, $\phi$, and $h$, such that there holds, for $0<\rho\leq \rho_0$,
\begin{equation}\label{eq:weyl}
\left| \rho^{n-1} \sum_{k' \in Q_{4/(3\rho)}'\cap\Z^{n-1}} \frac{\phi (\rho ( k',e'\cdot k'))}{\sum_{i=1}^{\hat{N}} \phi(  y^i + \rho (k',e'\cdot k'))} (h(e'\cdot k' + t) - \overline{h}) \right| \leq \eta,
\end{equation}
uniformly for any $t\in\T$, any $1\leq \hat{N}\leq N$ and any family $\{ y^i\}_{i=1}^{\hat{N}}\subset Q_2$ such that 
\begin{equation}\label{eq:xi-cond}
\sum_{i=1}^{\hat{N}} \phi( y^i + \rho (k',e' \cdot k')) \geq \beta \quad\text{for all }k'\in Q_{4/(3\rho)}'\cap\Z^{n-1}.
\end{equation}
\end{lemma}

\begin{proof} For simplicity, we will only deal with the case when $\hat{N}=N$. For notational convenience, let us write 
\begin{equation*}
\psi_\rho (k') = \psi_\rho (k'; e',\{ y^i\}_{i=1}^N) :=  \frac{\phi (\rho ( k',e'\cdot k'))}{\sum_{i=1}^N \phi(  y^i + \rho (k',e'\cdot k'))}.
\end{equation*}
Also denote by $I_\rho$ the set $Q_{4/(3\rho)}'\cap\Z^{n-1}$. 

\noindent{\it Step 1: The inequality \eqref{eq:weyl} holds for each $\{ y^i\}_{i=1}^N$.}

Let us consider a measure-preserving transformation group $\{T_{k'}\}_{k'\in\Z^{n-1}}$ on $\T$ defined by $T_{k'}(t) := t + e' \cdot k'$, and set 
\begin{equation}\label{eq:vprho}
\vp_\rho(t) := \rho^{n-1}\sum_{k'\in I_\rho} \psi_\rho(k') (h ( T_{k'} (t) ) - \overline{h})\quad (t\in\T).
\end{equation}
Note from the fact that $\overline{h \circ T_{k'} } = \overline{h}$, we have 
\begin{equation}\label{eq:vprho-avg}
\overline{\vp_\rho} = 0.
\end{equation}

We claim that 
\begin{equation}\label{eq:vprho-erg}
(\vp_\rho\circ T_{k'}  - \vp_\rho) \ra 0 \quad\text{uniformly on $\T$ for each $k'\in\Z^{n-1}$}.
\end{equation}
Fix $k_0'\in\Z^{n-1}$ and let us write by $E_\rho$ the symmetric difference $I_\rho \setminus (k_0' + I_\rho) \bigcup (k_0' + I_\rho) \setminus I_\rho$. By the associativity of $\{ T_{k'} \}_{k'\in\Z^{n-1}}$ and the triangle inequality, we have
\begin{equation*}
\begin{split}
|\vp_\rho \circ T_{k_0'}  - \vp_\rho| &\leq \rho^{n-1}\sum_{k' \in I_\rho}|\psi_\rho (k') - \psi_\rho (k'-k_0')||h \circ T_{k'+k_0'}  - \overline{h}| \\
& \quad + \rho^{n-1} \sum_{k' \in E_\rho} \psi_\rho(k') |h \circ T_{k'}  - \overline{h}|.
\end{split}
\end{equation*}
Let us take $\hat\rho>0$ small enough such that $2\alpha |k_0'||(e',1)|\hat\rho \leq \beta$. Then for all $0<\rho\leq \hat\rho$ and any $k\in I_\rho$, 
\begin{equation*}
\sum_{i=1}^N \phi ( y^i + \rho (k' - k_0' , e'\cdot (k'-k_0'))) \geq 2^{-1}\beta,
\end{equation*}
which yields that
\begin{equation*}
| \psi_\rho (k') - \psi_\rho (k'- k_0') |\leq 2\beta^{-2}\alpha |k_0'||(e',1)|\rho.
\end{equation*}
Thus, as we put $K:= \osc_\T h$, we obtain from the above inequalities and the fact that $\# E_\rho \leq |k_0'| \rho^{2-n}$, 
\begin{equation*}
|\vp_\rho \circ T_{k_0'} - \vp_\rho| \leq K|k_0'|(2\beta^{-2}\alpha|(e',1)| + 1)\rho\quad\text{for all }0<\rho\leq \hat\rho,
\end{equation*}
which leads us to \eqref{eq:vprho-erg}. It is noteworthy from the right-hand side of the last inequality that \eqref{eq:vprho-erg} does not depend on the choice of $\{ y^i\}_{i=1}^N$. 

Utilizing \eqref{eq:vprho-erg}, we are able to observe that $\{\vp_\rho\}_{\rho>0}$ converges uniformly to $0$ on $\T$ as $\rho\ra 0$. The equicontinuity and the uniform boundedness of $\{\vp_\rho\}_{\rho>0}$ are immediate from the definition \eqref{eq:vprho}, and the assumption on $h$ and $\phi$. Then we may deduce from Arzela-Ascoli theorem that for any subsequence $\{\vp_{N_i}\}_{i=1}^\infty$ of $\{\vp_\rho\}_{\rho>0}$, there exists a further subsequence $\{\vp_{M_j}\}_{i=1}^\infty\subset\{\vp_{N_i}\}_{i=1}^\infty$ and a function $\vp_0\in C(\T)$ such that $\vp_{M_j} \ra \vp_0$ uniformly on $\T$. Combining this observation with \eqref{eq:vprho-erg}, we deduce that 
\begin{equation}\label{eq:vp0-erg}
\vp_0\circ  T_{k'}  = \vp_0 \quad(k'\in\Z^{n-1}).  
\end{equation}
Now it follows from the ergodicity that $\vp_0$ is constant on $\T$. Then owing to \eqref{eq:vprho-avg}, we have
\begin{equation}\label{eq:vp0}
\vp_0 \equiv  \int_\T \vp_0 = \lim_{j\ra\infty} \int_\T \vp_{M_j} =0.
\end{equation}   

To summarize, we have shown that any subsequence of $\{\vp_\rho\}_{\rho>0}$ contains a further subsequence converging to $0$ uniformly on $\T$. In other words, $\vp_\rho\ra 0$ uniformly on $\T$ as $\rho \ra 0$. This implies that given $\eta>0$, there exists $\bar\rho$, which depends only on $\eta$, $n$, $\alpha$, $\beta$, $N$, $\{ y^i\}_{i=1}^N$, $\phi$ and $h$, such that \eqref{eq:weyl} holds for all $0<\rho\leq \bar\rho$. We are left to show that we may choose $\bar\rho$ independently on $\{ y^i\}_{i=1}^N$. 

\noindent{\it Step 2: Inequality \eqref{eq:weyl} holds uniformly for $\{ y^i\}_{i=1}^N$.}

Fix $\eta>0$ and take the least integer $m$ such that $m\geq\frac{1}{\eta}$ and $|h(t)-h(s)|\leq \eta$ for any $|t-s|\leq\frac{1}{m}$. Let $\cA_m$ be the family of all $L = \{l^i\}_{i=1}^N\subset Q_2\cap\frac{1}{m}\Z^n$ such that \eqref{eq:xi-cond} is satisfied for $ y^i$ replaced by $\{l^i\}_{i=1}^N$. We know that $\cA_\eta$ is a finite family whose cardinality is at most $m^{nN}$. 

Next for each $L\in\cA_m$, let us take $\rho_i$ as above such that \eqref{eq:weyl} holds with $\frac{\eta}{3}$ instead of $\eta$,  for any $0<\rho\leq\rho_i$. Collect them as a family $\{\rho_j\}_{j\in J}$. Then $\# J \leq m^{nN}$, whence we have that $\rho_0 := \min\{\rho_j:j\in J\} >0$. Note that $\rho_0$ depends only on $\eta$, $n$, $\alpha$, $\beta$, $N$, $\phi$ and $h$, and that \eqref{eq:weyl} holds for all $0<\rho\leq \rho_0$ and for all $L\in\cA_m$.

Now given any $X = \{ y^i\}_{i=1}^N\subset Q_2$ satisfying \eqref{eq:xi-cond}, we may choose some $\overline{X} = \{\overline{y}^i\}_{i=1}^N\in\cA_m$ such that $|y^i -  \overline{y}^i|\leq\frac{\sqrt{n}}{m}$ for each $1\leq i\leq N$, and $|\gamma-a|\leq \frac{1}{m}$. Then it is easy to observe that 
\begin{equation*}
\left| \frac{\phi (\rho ( k',e'\cdot k'))}{\sum_{i=1}^N \phi(  y^i + \rho (k',e'\cdot k'))} - \frac{\phi (\rho ( k',e'\cdot k'))}{\sum_{i=1}^N \phi( \overline{y}^i + \rho (k',e'\cdot k'))} \right| \leq \frac{\alpha N\sqrt{n}}{\beta^2m}. 
\end{equation*}
Thus, given $\eta>0$, taking $m$ large enough such that $m\geq \frac{K}{\alpha N\sqrt{n}3\beta^2\eta}$, we arrive at \eqref{eq:weyl} for any $0<\rho\leq \rho_0$, owing to the triangle inequality. This finishes the proof.
\end{proof}

The next lemma is an integral version of Lemma \ref{lemma:weyl}. 

\begin{lemma}\label{lemma:avg} Let $e'\in\R^{n-1}$ and $\phi\in C^1(\R^n)$ be as in Lemma \ref{lemma:weyl} and let $g\in C(\T^n)$ be given. Then for any $\eta>0$, there exists $\rho_0>0$, which depends only on $\eta$, $n$, $e'$, $\alpha$, $\beta$, $N$, $\phi$ and $g$, such that there holds, 
\begin{equation}\label{eq:avg}
\rho^{n-1}\left| \int_{Q_{4/(3\rho)}'} \frac{\phi (\rho ( y',e'\cdot y'))}{\sum_{i=1}^{\hat{N}} \phi(  y^i + \rho (y',e'\cdot y'))} (g((y',e'\cdot y') + y^0) - \overline{g}) dy' \right| \leq \eta,
\end{equation}
for $0<\rho\leq \rho_0$, uniformly for any $y^0\in\T^n$, any $\hat{N}\leq N$ and any family $\{ y^i\}_{i=1}^{\hat{N}}\subset Q_2$ satisfying \eqref{eq:xi-cond}.
\end{lemma}

\begin{proof} Let us define $h\in C(\T)$ by 
\begin{equation*}
h(t) := \int_{\T^{n-1}} g((y',e'\cdot y' + t) + y^0)dy' - \overline{g} \quad (t\in\T).
\end{equation*}
By the periodicity of $g$, $\overline{h} = 0$. We also observe that 
\begin{equation*}
\rho^{n-1}\left| \int_{Q_{4/(3\rho)}'} \psi_\rho (y') (g(y',e'\cdot y') - \overline{g}) dy' -  \sum_{k\in Q_{4/(3\rho)}'\cap\Z^{n-1}} \psi_\rho (k') h(e'\cdot k') dy' \right|\leq C\rho,
\end{equation*}
where $C$ depends only on $n$, $\alpha$, $\beta$, $N$ and $\osc_{\T^n} g$, and the function $\psi_\rho: Q_{4/(3\rho)}'\ra \R$ is defined by $\psi_\rho(y') := \phi(y',e'\cdot y')/ \sum_{i=1}^{\hat{N}} \phi( y^i + \rho (y',e'\cdot y'))$. The proof is now finished by Lemma \ref{lemma:weyl}.
\end{proof}

\subsection{Averaging Behavior under IDDC}\label{subsection:proof-osc}

The main objective of this subsection is to prove Theorem \ref{theorem:osc-int}. As the first step, let us consider the case when $\Gamma\in C^1$ and $g\in C(\Gamma\times\T^n)$. 

\begin{proposition}\label{proposition:osc-int} Let $\Gamma\subset\R^n$ ($n\geq 2$) be a compact $C^1$-hypersurface and $g\in C(\Gamma\times\T^n)$. Then one has the convergence in \eqref{eq:osc-int}, provided $\Gamma$ satisfies IDDC. 
\end{proposition}

Before going into the proof, let us introduce some notation which we will use in the proof of the above proposition. Due to the $C^1$-regularity of $\Gamma$, $\nu$ is continuous on $\Gamma$, and moreover, we may use the compactness of $\Gamma$ to quantify the continuity of $\nu$ by invoking a modulus of continuity $\tau$; that is,
\begin{equation}\label{eq:Gamma-C1}
| \nu(z) - \nu(\overline{z}) | \leq \tau ( |z-\overline{z}| )\quad\text{for any }z,\overline{z} \in\Gamma.
\end{equation}
Recall from \eqref{eq:Pi} that $\Pi(z)$ is the tangent hyperplane at $z$. We may choose $r_0>0$ such that the projection mapping $\pi_{z,r}$ from $\Pi(z)   \cap Q_{4r/3}(z)$ onto $\Gamma\cap Q_{4r/3}(z)$, defined by 
\begin{equation}\label{eq:pizr}
\pi_{z,r}(x) = x + \dist(x,\Gamma) \nu(z),
\end{equation}
is a $C^1$-diffeomorphism, for all $0<r\leq r_0$. By \eqref{eq:Gamma-C1}, it is not hard to observe that 
\begin{equation}\label{eq:pizr-C1}
| x - \pi_{z,r}(x) | \leq cr\tau\left(\frac{4r}{3}\right)\quad\text{and}\quad | I - \nabla \pi_{z,r} | \leq c\tau\left(\frac{4r}{3}\right),
\end{equation}
where $c$ is a dimensional constant, for any $0<r\leq r_0$. Let us note that a similar argument can be found also in \cite{AKMP} and \cite{GM2}.

Without loss of generality, we may assume that $\tau$ controls the continuity of $g$ as well; i.e., 
\begin{equation}\label{eq:g-C}
|g(z,y) - g(\overline{z},\overline{y})| \leq \tau( (|z-\overline{z}|^2 + |y-\overline{y}|^2)^{1/2})\quad\text{for any }(z,y),(\overline{z},\overline{y})\in\Gamma\times\T^n.
\end{equation}
Let us define, for notational simplicity, $\overline{g}:\Gamma\ra \R$ by 
\begin{equation}\label{eq:gov-osc-int}
\overline{g}(z) := \int_{\T^n} g(z,y)dy.
\end{equation}
It is evident from \eqref{eq:g-C} that $\overline{g}\in C(\Gamma)$ and 
\begin{equation}\label{eq:gov-C1}
|\overline{g}(z) - \overline{g}(\overline{z})| \leq \tau(|z-\overline{z}|)\quad\text{for any }z,\overline{z}\in\Gamma.
\end{equation}
Let us take a positive number $K$ such that
\begin{equation}\label{eq:K}
K \geq \max \left\{ \sigma(\Gamma),\norm{g}_{L^\infty(\Gamma\times\T^n)},1 \right\} \geq \max\left\{\sigma(\Gamma),\norm{\overline{g}}_{L^\infty(\Gamma)}, 1\right\}
\end{equation}

We are now in position to prove Proposition \ref{proposition:osc-int}.

\begin{proof}[Proof of Proposition \ref{proposition:osc-int}]
Throughout this proof, we will use the notations $\nu$, $\Pi(z)  $, $\pi_{z,r}$, $\overline{g}$ and $K$ introduced in \eqref{eq:Gamma-C1} -- \eqref{eq:K}. Let us fix $\eta$ as an arbitrary positive number, $N$ a large dimensional constant and $\phi\in C^1(\R^n)$ a cut-off function satisfying
\begin{equation}\label{eq:phi}
\begin{cases}
0\leq \phi\leq 1\text{ and }|\nabla \phi|\leq C\text{ on }\R^n,\\
\phi\equiv 1\text{ on }Q_1\text{ and }\spt\phi\subset Q_{\frac{4}{3}},
\end{cases}
\end{equation}
where $C$ is a constant depending only on $n$. In spite of slight abuse of notation, we assume for simplicity that $\nu_n(z) \neq 0$ for all $z\in\Gamma$, and define $e':\Gamma\ra\R^{n-1}$ by 
\begin{equation}\label{eq:e'}
e'(z) := -\nu_n(z)^{-1} \nu'(z).
\end{equation} 

\noindent{\it Step 1: Determination of ratio $\rho$ of microscopic scale $\e$ to mesoscopic scale $r$.} 

Given $\phi$, let us consider a measurable function $F$ on $\Gamma\times(0,1)$ defined by 
\begin{equation}\label{eq:F}
F(z,\rho):= \fint_{Q_{\frac{4}{3\rho}}'} \frac{ \phi( \rho (y',e'(z)\cdot y')) }{\sum_{i=1}^{N'} \phi(  y^i + \rho(y',e'(z)\cdot y') )} ( g(z,(y',e'(z)\cdot y') + y^0) - \overline{g}(z) ) dy',
\end{equation}
where $y^0\in\T^n$, $N'\leq N$ and $\{ y^i\}_{i=1}^{N'}\subset Q_2$ such that \eqref{eq:xi-cond} is true with $\alpha=C$, $\beta=\frac{1}{2}$ and $\gamma =2N$ and $N$ replaced by $\frac{1}{2}$ and respectively $2N$. 

By definition, $\nu(z) = \frac{(e'(z),1)}{|(e'(z),1)|}\in\Ss^{n-1}$ for any $z\in\Gamma$. Thus, it follows from IDDC that $F(z,\rho )\ra 0$ as $\rho\ra 0$ for a.e. $z\in\Gamma$; more specifically, we know from \eqref{eq:avg} that given $\eta>0$, for a.e. $z\in\Gamma$, there exists $\rho_{z,\eta}$, depending only on $\eta$, $n$, $z$, $\phi$ and $g$, such that
\begin{equation}\label{eq:F-limit}
| F(z,\rho) | \leq \eta\quad\text{for all }0<\rho\leq \rho_z.
\end{equation}
Note that $\rho_z$ is chosen uniform over any $y^0\in\T^n$, $N'\leq N$ and $\{ y^i\}_{i=1}^{N'}$ given as above. Hence, utilizing Egoroff's theorem, we obtain, for each $\eta>0$, a subset $\Gamma_0\subset \Gamma$ such that 
\begin{equation}\label{eq:offset}
\text{ $\sigma(\Gamma\setminus\Gamma_0)\leq \eta$ and $F(z,\rho)\ra 0$ as $\rho\ra 0$ uniformly for $z\in\Gamma_0$,}
\end{equation}
and accordingly choose a number $\rho_0>0$ for which 
\begin{equation}\label{eq:meso}
\norm{F(\cdot,\rho)}_{L^\infty(\Gamma_0)}\leq \eta\quad\text{for all $0<\rho\leq \rho_0$}.
\end{equation}

\noindent{\it Step 2: Construction of a Vitali covering and a partition of unity of $\Gamma_0$.} 

Let us fix $r>0$ and consider an open covering $\{Q_r(z)\}_{z\in\Gamma_0}$ of $\Gamma_0$. By the Vitali covering lemma, there is a countable set $\cA_r\subset\Gamma_0$ such that 
\begin{equation}\label{eq:vitali-disjoint}
\Gamma_0 \subset \bigcup_{z\in\cA_r} Q_r(z)\quad\text{and}\quad Q_{\frac{r}{3\sqrt{n}}}(z)\cap Q_{\frac{r}{3\sqrt{n}}}(\overline{z})= \emptyset\quad\text{for any $z,\overline{z}\in\cA_r$ with $z\neq \overline{z}$}.
\end{equation}
We may refine $\cA_r$ to satisfy 
\begin{equation}\label{eq:vitali-intersect}
1\leq \#\{ \overline{z}\in\cA_r\setminus\{z\}: Q_{4r/3}(z) \cap Q_r(\overline{z})\neq\emptyset\} \leq N,
\end{equation}
by taking $N$ larger if necessary; indeed we choose the dimensional constant $N$ to satisfy \eqref{eq:vitali-intersect}. For notational convenience, let us denote 
\begin{equation}\label{eq:Nr}
\cN_r := \bigcup_{z\in\cA_r} Q_r(z).
\end{equation} 

Next we define 
\begin{equation}\label{eq:Phir}
\Phi_r(x) := \sum_{z\in\cA_r} \phi \left( \frac{ x- z }{r} \right) \quad\text{on }\R^n.
\end{equation}
Due to \eqref{eq:vitali-intersect}, the summation on the right-hand side involves at most $N$-terms. It also follows from \eqref{eq:phi}, \eqref{eq:vitali-disjoint} and \eqref{eq:vitali-intersect} that 
\begin{equation}\label{eq:Phir-C1}
1\leq \Phi_r \leq N\quad\text{on }\cN_r \quad\text{and}\quad r|\nabla \Phi_r|\leq CN\quad\text{on }\R^n. 
\end{equation}
In what follows, we define, for each $z\in\cA_r$, 
\begin{equation}\label{eq:psizr}
\psi_{z,r}(x) := \frac{1}{\Phi_r(x)}\phi\left(\dfrac{x-z}{r}\right)\quad\text{on }\cN_r.
\end{equation} 
Owing to \eqref{eq:phi} and \eqref{eq:Phir}, we observe that $\psi_{z,r}$ satisfies 
\begin{equation}\label{eq:partition}
\begin{dcases}
0\leq \psi_{z,r}\leq 1\text{ and } r |\nabla \psi_{z,r}| \leq CN\text{ on }\cN_r,\\
\spt(\psi_{z,r}) \subset Q_{\frac{4r}{3}}(z),\\
\sum_{z\in\cA_r} \psi_{z,r} =1\quad\text{on }\cN_r.
\end{dcases}
\end{equation}

\noindent{\it Step 3: Decomposition of the boundary $\Gamma$.} 

Let us take $r_0>0$, chosen in Step 1, smaller if necessary such that $cr\tau(4r/3)\leq\frac{1}{2}$ for all $0<r\leq r_0$. Then it follows from \eqref{eq:pizr-C1} and \eqref{eq:Phir-C1} that
\begin{equation}\label{eq:Phir-Pi}
\frac{1}{2}\leq \Phi_r \leq 2N\quad\text{on } \Pi(z)   \cap Q_{\frac{4r}{3}}(z),
\end{equation}
uniformly for all $z\in\Gamma$. On the other hand, owing to \eqref{eq:vitali-intersect}, there is $N_{z,r}\leq N$ for each $z\in\cA_r$ such that $\{ \overline{z}\in\cA_r\setminus\{z\}: Q_{4r/3}(z) \cap Q_r(\overline{z})\neq\emptyset\}=\{\overline{z}_r^1,\overline{z}_r^2,\cdots, \overline{z}_r^{N_{z,r}}\}$. Let us denote, for each $1\leq i\leq N_{z,r}$, 
\begin{equation}\label{eq:xizr}
y_{z,r}^i := \frac{z - \overline{z}_r^i}{r} \in Q_2,
\end{equation}
where the inclusion follows from \eqref{eq:vitali-intersect} and is true for any $z\in\cA_r$ and any $r>0$. 

Now we are in a position to introduce the microscopic scale parameter $\e$ and put $r= \frac{\e}{\rho_0}$, where $\rho_0$ is chosen from \eqref{eq:meso}. Then keeping in mind the transformation $x = z + \e (y', e'(z)\cdot y')$ from $y'\in Q_{\frac{4}{3\rho_0}}'$ onto $ \Pi(z)   \cap Q_{\frac{4\e}{3\rho_0}}(z)$, we may rephrase \eqref{eq:Phir-Pi} as 
\begin{equation}\label{eq:Phir-Pi-re}
\frac{1}{2} \leq \Phi_{\frac{\e}{\rho_0}} (x) = \sum_{i=1}^{N_{z,\e/\rho_0}} \phi (y_{z,\frac{\e}{\rho_0}}^i + \rho_0 (y',e'(z)\cdot y') ) \leq 2N,
\end{equation}
for all $\e>0$. Thus, by substituting $\frac{z}{\e}$ and $\{y_{z,\e/\rho_0}^i\}_{i=1}^{N_{z,\e/\rho_0}}$ into $y^0$ and respectively $\{ y^i\}_{i=1}^{N'}$ in the definition \eqref{eq:F} of $F$, it follows from \eqref{eq:meso} that 
\begin{equation}\label{eq:meso-re}
\left| F\left(z,\rho_0; \frac{z}{\e},\{y_{z,\frac{\e}{\rho_0}}^i\}_{i=1}^{N_{z,\frac{\e}{\rho_0}}}\right) \right| \leq \eta\quad\text{for any $z\in\cA_{\frac{\e}{\rho_0}}$ and any $\e>0$}.
\end{equation}
Utilizing \eqref{eq:Phir-Pi-re}, we derive that 
\begin{equation}\label{eq:F-re}
F\left(z,\rho_0; \frac{z}{\e},\{y_{z,\frac{\e}{\rho_0}}^i\}_{i=1}^{N_{z,\frac{\e}{\rho_0}}}\right) = \fint_{ \Pi(z)   \cap Q_{\frac{4\e}{3\rho_0}}(z)} \psi_{z,\frac{\e}{\rho_0}}(x) \left( g\left(z,\frac{x}{\e}\right) - \overline{g}(z)\right)d\sigma_x,
\end{equation}
for any $z\in\cA_{\frac{\e}{\rho_0}}$ and any $\e>0$.

\noindent{\it Step 4: Decomposition of the oscillatory integral $\int_\Gamma g(x,\e^{-1}x)d\sigma_x$.}

Let $\mu$ be a universal constant to be determined and set $\e_0\leq \min\{ \rho_0 r_0,1\}$ such that 
\begin{equation}\label{eq:e0}
\max\left\{ \tau\left(\frac{4\e}{3\rho_0}\right) , \tau \left( \frac{\mu}{\rho_0} \tau \left(\frac{4\e}{3\rho_0} \right)\right) \right\} \leq \eta\quad\text{for any }0<\e\leq \e_0.
\end{equation}
Let us decompose the integral $\int_\Gamma(g(x,\frac{x}{\e}) - \overline{g}(x))d\sigma_x$ as follows:
\begin{equation}\label{eq:decom}
\left| \int_\Gamma g\left(x,\frac{x}{\e}\right) d\sigma_x - \int_\Gamma \overline{g}(x) d\sigma_x \right| \leq I_1 + I_2 + I_3 + I_4 + I_5,
\end{equation}
where
\begin{align}
\label{eq:I1}
I_1 &:= \left| \int_\Gamma g\left(x,\frac{x}{\e}\right) d\sigma_x - \sum_{z\in \cA_{\e/\rho_0}} \int_{\Gamma\cap Q_{\frac{4\e}{3\rho_0}}(z)} \psi_{z,\frac{\e}{\rho_0}} (x) g\left(x,\frac{x}{\e}\right) d\sigma_x \right|, \\
\label{eq:I2}
I_2 &:= \sum_{z\in\cA_{\e/\rho_0}} \left| \int_{\Gamma\cap Q_{\frac{4\e}{3\rho_0}}(z)} \psi_{z,\frac{\e}{\rho_0}}(x) g\left(x,\frac{x}{\e}\right) d\sigma_x - \int_{ \Pi(z)   \cap Q_{\frac{4\e}{3\rho_0}}(z)} \psi_{z,\frac{\e}{\rho_0}} (x) g\left(z,\frac{x}{\e}\right) d\sigma_x \right|, \\
\label{eq:I3}
I_3 &:= \sum_{z\in\cA_{\e/\rho_0}} \left| \int_{ \Pi(z)   \cap Q_{\frac{4\e}{3\rho_0}}(z)} \psi_{z,\frac{\e}{\rho_0}}(x) \left( g\left(z,\frac{x}{\e}\right) - \overline{g}(z)\right)d\sigma_x \right|,\\
\label{eq:I4}
I_4 &:= \sum_{z\in\cA_{\e/\rho_0}} \left| \overline{g}(z) \int_{ \Pi(z)   \cap Q_{\frac{4\e}{3\rho_0}}(z)} \psi_{z,\frac{\e}{\rho_0}}(x) d\sigma_x - \int_{\Gamma\cap Q_{\frac{4\e}{3\rho_0}}(z)} \psi_{z,\frac{\e}{\rho_0}}(x)\overline{g}(x) d\sigma_x \right|,
\end{align}
and 
\begin{equation}\label{eq:I5}
I_5 := \left| \sum_{z\in\cA_{\e/\rho_0}} \int_{\Gamma\cap Q_{\frac{4\e}{3\rho_0}}(z)} \psi_{z,\frac{\e}{\rho_0}}(x)\overline{g}(x) d\sigma_x - \int_\Gamma \overline{g}(x) d\sigma_x\right|.
\end{equation}

\noindent{\it Step 5: Estimation of $I_1$ and $I_5$.} 

One should notice that the partition of unity in Step 2 is performed in the set $\Gamma_0\subset\Gamma$ defined in \eqref{eq:offset}.  Since
\begin{equation*}
\int_{\Gamma_0} g\left(x,\frac{x}{\e}\right) d\sigma_x = \sum_{z\in\cA_{\e/\rho_0}} \int_{\Gamma_0\cap Q_{\frac{4\e}{3\rho_0}}(z)} \psi_{z,\frac{\e}{\rho_0}} (x) g\left(x,\frac{x}{\e}\right) d\sigma_x,
\end{equation*}
we may deduce from the first inequality in \eqref{eq:offset} that
\begin{equation}\label{eq:decom1} 
I_1 \leq 2K\eta. 
\end{equation}
In the similar way, we deduce that 
\begin{equation}\label{eq:decom5}
I_5 \leq 2K\eta. 
\end{equation}

\noindent{\it Step 6: Estimation of $I_2$ and $I_4$.}

This part mainly involves the regularity of $\Gamma$, $g$ and $\overline{g}$, but has nothing to do with ergodicity. Let us fix $z\in\cA_{\e/\rho_0}$. As we write by $J_{z,\frac{\e}{\rho_0}}$ the Jacobian determinant of $\pi_{z,\frac{\e}{\rho_0}}^{-1}$, we have 
\begin{equation*}
\begin{split}
& \int_{\Gamma\cap Q_{\frac{4\e}{3\rho_0}}(z)} \psi_{z,\frac{\e}{\rho_0}} (x) g\left( x, \frac{x}{\e} \right) d\sigma_x \\
& = \int_{\Pi(z)  \cap Q_{\frac{4\e}{3\rho_0}}(z)} \psi_{z,\frac{\e}{\rho_0}} (\pi_{z,\frac{\e}{\rho_0}}(x)) g \left( \pi_{z,\frac{\e}{\rho_0}}(x), \frac{\pi_{z,\frac{\e}{\rho_0}}(x)}{\e} \right) J_{z,\frac{\e}{\rho_0}}(x) d\sigma_x.
\end{split}
\end{equation*}
Therefore, the estimation of $I_2$ essentially follows from $\psi_{z,\frac{\e}{\rho_0}}(x) - \psi_{z,\frac{\e}{\rho_0}}(\pi_{z,\frac{\e}{\rho_0}}(x))$, $g(x,\frac{x}{\e}) - g ( \pi_{z,\frac{\e}{\rho_0}}(x), \e^{-1}\pi_{z,\frac{\e}{\rho_0}}(x) )$ and $J_{z,\frac{\e}{\rho_0}}(x) - 1$ on $\Pi(z)  \cap Q_{\frac{4\e}{3\rho_0}}(z)$. It follows from \eqref{eq:pizr-C1} and \eqref{eq:partition} that
\begin{equation*}
\left| \psi_{z,\frac{\e}{\rho_0}}(x) - \psi_{z,\frac{\e}{\rho_0}}(\pi_{z,\frac{\e}{\rho_0}}(x)) \right| \leq c \norm{\nabla \psi_{z,\frac{\e}{\rho_0}}}_{L^\infty(\Pi(z)  \cap Q_{\frac{4\e}{3\rho_0}}(z))}\left(\frac{\e}{\rho_0}\right)\tau\left( \frac{4\e}{3\rho_0} \right) \leq C\tau\left(\frac{4\e}{3\rho_0}\right).
\end{equation*}
Similarly, we utilize \eqref{eq:pizr-C1} and \eqref{eq:g-C} to derive that
\begin{equation*}
\left| g\left(x,\frac{x}{\e}\right) - g \left( \pi_{z,\frac{\e}{\rho_0}}(x), \frac{\pi_{z,\frac{\e}{\rho_0}}(x)}{\e} \right) \right| \leq \tau\left( \left(1 + \frac{1}{\e^2}\right)^{1/2} \left|x - \pi_{z,\frac{\e}{\rho_0}}(x)\right| \right) \leq \tau \left(\frac{\mu}{\rho_0} \tau\left(\frac{4\e}{3\rho_0}\right) \right),
\end{equation*} 
provided $\e\leq 1$ and $\mu\geq c\sqrt{2}$ with $c$ chosen from \eqref{eq:pizr-C1}. On the other hand, it immediately follows from \eqref{eq:pizr-C1} that 
\begin{equation*}
\left| J_{z,\frac{\e}{\rho_0}}(x) - 1 \right| \leq C\tau\left(\frac{4\e}{3\rho_0}\right). 
\end{equation*}
In particular, we have $|J_{z,\frac{\e}{\rho_0}}| \leq 2$ by choosing $\e$ so small that $C\tau(\frac{4\e}{3\rho_0})\leq 1$. Collecting the three inequalities above and utilizing \eqref{eq:vitali-intersect}, we arrive at 
\begin{equation}\label{eq:decom2}
I_2 \leq (3C +2)KN\eta,  
\end{equation}
due to the choice of $\e_0$ made in \eqref{eq:e0}.

One may obtain a similar estimate for $I_4$; indeed, as $\overline{g}$ being independent of variable $y$, the estimate becomes simpler. We observe that
\begin{equation}\label{eq:decom4}
I_4 \leq (3C+2)KN\eta,
\end{equation}
for which we skip the details. 

\noindent{\it Step 7: Estimation of $I_3$.}

Owing to \eqref{eq:meso-re} and \eqref{eq:F-re}, we know that for any $z\in\cA_{\e/\rho_0}$ and any $0<\e\leq \e_0$, there holds 
\begin{equation*}
\left| \fint_{ \Pi(z)   \cap Q_{\frac{4\e}{3\rho_0}}(z)} \psi_{z,\frac{\e}{\rho_0}}(x) \left( g\left(z,\frac{x}{\e}\right) - \overline{g}(z)\right)d\sigma_x \right| \leq \eta. 
\end{equation*}
Therefore,
\begin{equation}\label{eq:decom3}
I_3 \leq \left[ \sum_{z\in\cA_{\e/\rho_0}} \sigma\left(\Pi(z)  \cap Q_{\frac{4\e}{3\rho_0}}(z)\right) \right] \eta \leq CNK\eta,
\end{equation}
since $\sigma(\Pi(z)  \cap Q_{\frac{4\e}{3\rho_0}} (z))\leq 2\sigma(\Gamma\cap Q_{\frac{4\e}{3\rho_0}})$, which follows from $|J_{z,\frac{\e}{\rho_0}}| \leq 2$, and since the elements of the covering $\{Q_{\frac{4\e}{3\rho_0}}(z)\}_{z\in\cA_{\frac{\e}{\rho_0}}}$ mutually intersect at most $CN$-times, which is ensured by \eqref{eq:vitali-intersect}.

Finally, collecting inequalities \eqref{eq:decom1}, \eqref{eq:decom2}, \eqref{eq:decom3}, \eqref{eq:decom4} and \eqref{eq:decom5} altogether, we deduce from \eqref{eq:vitali-intersect} and \eqref{eq:e0} that
\begin{equation}\label{eq:decom}
\left| \int_\Gamma g\left( x, \frac{x}{\e} \right) d\sigma_x - \int_\Gamma \overline{g}(x) d\sigma_x \right| \leq ((7C+4)N+4)K\eta, 
\end{equation}
provided that $0<\e\leq \e_0$, finishing the proof of Proposition \ref{proposition:osc-int}.
\end{proof}

We are now ready to prove Theorem \ref{theorem:osc-int}. With Proposition \ref{proposition:osc-int} at hand, we make a careful approximation argument for $g$, and apply the preceding proposition to have the desired limit. 

\begin{proof}[Proof of Theorem \ref{theorem:osc-int}]
Denote by $\tau$ the modulus of continuity of $g(x,\cdot)$ on $\T^n$. We claim that there exists a compact $\Gamma_0\subset\Gamma$ for which 
\begin{equation}\label{eq:claim-osc-int}
\text{$\sigma(\Gamma\setminus\Gamma_0)<\eta$ and $g$ is continuous on $\Gamma_0\times\T^n$.}
\end{equation}

Let $k$ be a positive integer and split $\T^n$ into open cubes $Q_j$, for $j=1,2,\cdots,2^{kn}$, where $Q_j$ has side length $2^{-k}$ and $\bigcup_{j=1}^{2^{kn}} \overline{Q}_j = \T^n$. Next we fix $y^j\in Q_j$ and by Lusin's theorem, choose a closed $\Gamma_j\subset\Gamma$ such that $\sigma(\Gamma\setminus\Gamma_j)<2^{-kn}\eta$ and $g(\cdot,y^j)$ is continuous on $\Gamma_j$. 

Denote by $\tau_j$ the modulus of continuity for $g(\cdot,y^j)$ on $\Gamma_j$ and by $\tau$ the modulus of continuity for $g(x,\cdot)$ on $\T^n$; note that the equicontinuity of the family $\{g(x,\cdot)\}_{x\in\Gamma}$ allows us to choose $\tau$ uniformly over $x\in\Gamma$. Given $x^i\in\Gamma_j$ for $i=1,2$, we have $|g(x^i,y) - g(x^i,y^j)| \leq \tau(n^{1/2}2^{-k})$ and $|g(x^1,y^j) - g(x^2,y^j)|\leq \tau_j(|x^1-x^2|)$. Now taking $\tilde\tau = \max\{\tau_j:1\leq j\leq 2^{kn}\}$ and $|x^1-x^2|<\delta$, we observe from the triangle inequality that
\begin{equation*}
|g(x^1,y) - g(x^2,y)| \leq 2\tau(n^{1/2}2^{-k}) + \tilde\tau(\delta).
\end{equation*} 
As $\tilde\tau$ being a modulus of continuity, we are able to take $\delta$ small enough such that the rightmost side above becomes less than $3\tau(n^{1/2}2^{-k})$. Since $k$ is chosen to be an arbitrary positive integer, we deduce that $g(\cdot,y)$ is also continuous on $\Gamma_j$. 

Define $\Gamma_0 = \bigcap_{j=1}^{2^{kn}} \Gamma_j$. Since $\Gamma_j$ is closed, $\Gamma_0$ is also closed, and since $\Gamma$ is compact and $\Gamma_0\subset\Gamma$, we obtain that $\Gamma_0$ is compact. Moreover, it follows from $\sigma(\Gamma\setminus \Gamma_j)<2^{-kn}\eta$ for all $1\leq j\leq 2^{kn}$ that $\sigma(\Gamma\setminus \Gamma_0) < \eta$. Furthermore, as $g(\cdot,y)$ being continuous on $\Gamma_j$ for any $y\in \overline{Q}_j$, we have that $g(\cdot,y)$ is continuous on $\Gamma_0$ for any $y\in\T^n$. This proves \eqref{eq:claim-osc-int}.  

Let $\tilde{g}$ be a continuous extension of $g$ to $\Gamma\times\T^n$ such that 
\begin{equation}\label{eq:extension}
\norm{\tilde{g}}_{L^\infty(\Gamma\times\T^n)} \leq K := \norm{g}_{L^\infty(\Gamma\times\T^n)}.
\end{equation}
Proposition \ref{proposition:osc-int} is applicable to $\tilde{g}$ from which we obtain $\e_0>0$ such that
\begin{equation}\label{eq:gt-osc-int}
\left| \int_\Gamma \tilde{g} \left(x,\frac{x}{\e}\right) d\sigma_x - \int_\Gamma \int_{\T^n} \tilde{g}(x,y) dyd\sigma_x \right| \leq \eta,
\end{equation}
for all $0<\e\leq \e_0$. Thus, it follows from \eqref{eq:claim-osc-int}, \eqref{eq:extension} and the triangle inequality that
\begin{equation*}
\left| \int_\Gamma g \left(x,\frac{x}{\e}\right) d\sigma_x - \int_\Gamma \int_{\T^n} g(x,y) dyd\sigma_x \right| \leq (4K+1)\eta,
\end{equation*}
proving the theorem.
\end{proof}

%%%%%%%%%%%%%%%%%%%%%%%%%%%%%%%%%%%%%%%%%%%%%%%%%%%%%%%%%%%%%%%%%%%%%%%%%%%%%%%%%%%%%%%%%%%%%%%%
%
%     
%
%%%%%%%%%%%%%%%%%%%%%%%%%%%%%%%%%%%%%%%%%%%%%%%%%%%%%%%%%%%%%%%%%%%%%%%%%%%%%%%%%%%%%%%%%%%%%%%%

\section{Application to Homogenization of Oscillatory Boundary Data}\label{section:bdry-hom}

This section is devoted to the study of homogenization of \eqref{eq:Le}, as an application of the averaging behavior of oscillatory surface integrals established in Section \ref{section:osc-int}. We shall note that the main ideas of this section are scattered in different papers such as \cite{AL1}, \cite{AL2}, \cite{GM1}, \cite{GM2} and \cite{P}, and some of them central to our results can be found in \cite{AKMP}. 

Let us first make specific the assumptions in Theorem \ref{theorem:bdry-hom}. Let $A = (A_{ij}^{\alpha\beta})_{1\leq i,j\leq n}^{1\leq \alpha,\beta\leq m} : \T^n\ra \R^{n^2\times m^2}$ be uniformly elliptic in the sense of \eqref{eq:ellip-A} and $C^{ 1 + [\frac{n}{2}]}$-regular (see the definition of $[\cdot]$ in Section \ref{section:notation}) with
\begin{equation}\label{eq:Cka-A}
\norm{A}_{C^{ 1+ [\frac{n}{2}]}(\T^n;\R^{n^2\times m^2})} \leq K.
\end{equation}
Here $n\geq 2$ is the spatial dimension and $m\geq 1$ is the number of equations involved in \eqref{eq:Le}. Let us also assume that $\Omega\subset\R^n$ is a bounded $C^{1,\mu}$-domain whose boundary $\Gamma$ satisfies IDDC and 
\begin{equation}\label{eq:C1a-bdry}
|\nu(z) - \nu(\overline{z})| \leq K|z-\overline{z}|^\mu\quad\text{for all }z,\overline{z}\in\Gamma,
\end{equation}
where $\nu:\Gamma\ra \Ss^{n-1}$ is the unit inward normal mapping on $\Gamma$. Moreover, suppose that $g:\Gamma\times\T^n\ra \R^m$ is a bounded measurable mapping such that 
\begin{equation}\label{eq:Linf-g}
\norm{g}_{L^\infty(\Gamma\times\T^n;\R^m)}\leq K,
\end{equation}
and there is a modulus of continuity $\tau$ such that 
\begin{equation}\label{eq:equi-g}
\sup_{z\in\Gamma} |g(z,y) - g(z,\overline{y})| \leq \tau(|y-\overline{y}|)\quad\text{for any }y,\overline{y}\in\T^n. 
\end{equation}

Throughout this section, $c$ and $C$ are used to denote positive constants depending at most on $n$, $\lambda$, $K$, $\diam(\Omega)$, $\mu$ and $\tau$, where $\lambda$ is the ellipticity constant appearing in \eqref{eq:ellip-A}. Moreover, we will denote by $I$ the identity matrix in $\R^{m^2}$ and by $A^*$ the adjoint of $A$. Note that the adjoint $A^*$ is given by $(A^*)_{ij}^{\alpha\beta} = A_{ji}^{\beta\alpha}$ for any $1\leq i,j\leq n$ and $1\leq \alpha,\beta\leq m$, and that it satisfies \eqref{eq:ellip-A} and \eqref{eq:Cka-A} as well. Furthermore, we will fix, throughout this section, $\gamma$ as a dimensional constant given by
\begin{equation}\label{eq:gamma}
\gamma = 
\begin{dcases}
\left[\frac{n}{2}\right] +1 - \frac{n}{2},&\text{if $n$ is odd},\\
\text{any positive number $<1$},&\text{if $n$ is even}.
\end{dcases}
\end{equation}

\subsection{Constant Elliptic Systems}\label{subsection:example}

Let us begin with the elliptic system with constant coefficients, 
\begin{equation}\tag{$\overline{L}_\e$}\label{eq:Love}
\begin{dcases}
- \nabla \cdot (\overline{A} \nabla \overline{u}_\e(x)) = 0 & \text{in }\Omega, \\
\overline{u}_\e(x) = g\left(x,\frac{x}{\e}\right) & \text{on }\Gamma,
\end{dcases}
\end{equation}
where $\overline{A} \in \R^{n^2\times m^2}$ satisfies ellipticity condition \eqref{eq:ellip-A}. We observe that the boundary layer homogenization for \eqref{eq:Love} is a direct consequence of Theorem \ref{theorem:osc-int}, via the Poisson integral representation of $\overline{u}_\e$. Due to the absence of interior oscillation effect near the boundary layer, the effective boundary data $\overline{g}:\Gamma\ra\R^m$ is given by 
\begin{equation}\label{eq:gov-const}
\overline{g}(z) := \int_{\T^n} g(z,y) dy,
\end{equation}
which is only a simple average of $g$ in its second variable; compare this with \eqref{eq:gov} when interior oscillation is present. Additionally, the regularity assumption on $\Gamma$ can be relieved to continuous differentiability, rather than $C^{1,\mu}$-regularity. 

\begin{proposition}\label{proposition:const} Let $\overline{A}\in\R^{n^2\times m^2}$ be an elliptic coefficient in the sense of \eqref{eq:ellip-A}, $\Omega\subset\R^n$ a bounded $C^1$-domain whose boundary $\Gamma$ satisfies IDDC, and $g\in L^\infty(\Gamma\times\T^n;\R^m)$ satisfies that $\{g(z,\cdot)\}_{z\in\Gamma}$ is an equicontinuous family on $\T^n$. Then the sequence $\{\overline{u}_\e\}_{\e>0}$ of weak solutions of \eqref{eq:Love} converges, locally uniformly in $\Omega$, to the weak solution $\overline{u}$ of \eqref{eq:Lov} with $\overline{g}$ given by \eqref{eq:gov-const}. 
\end{proposition}

\begin{proof}
Let us denote by $\overline{P}$ the Poisson kernel for $-\nabla \cdot (\overline{A}\nabla )$ in $\Omega$, so that we have
\begin{equation*}
\overline{u}_\e(x) = \int_\Gamma \overline{P}(x,z) g\left(z,\frac{z}{\e}\right) d\sigma_z\quad(x\in\Omega).
\end{equation*}
Now fix $E\Subset \Omega$. Then the mapping $(z,y)\mapsto \overline{P}(x,z)g(z,y)$ on $\Gamma\times\T^n$ is bounded measurable whose bound is uniform over $E$, and that $\{\overline{P}(x,z)g(z,\cdot)\}_{z\in\Gamma,x\in E}$ is equicontinuous on $\T^n$. Owing to IDDC assumption on $\Gamma$, Theorem \ref{theorem:osc-int} is applicable to each component of the mapping $(z,y)\mapsto \overline{P}(x,z)g(z,y)$ on $\Gamma\times\T^n$; here we should remark that $\overline{P}(x,z)g(z,y) = (\overline{P}^{\alpha\beta} (x,z) g^\beta(x,z))_{1\leq \alpha\leq m}$. Thus,  
\begin{equation*}
\lim_{\e\ra 0} \overline{u}_\e(x) =  \int_\Gamma\int_{\T^n} \overline{P}(x,z)g(z,y)dyd\sigma_z = \int_\Gamma \overline{P}(x,z)\overline{g}(z) d\sigma_z = \overline{u}(x).
\end{equation*}
It should also be noted that the convergence above is uniform over all $x\in E$, since the convergence \eqref{eq:osc-int} relies only on the continuity and the supremum bound for $(z,y)\mapsto P(x,z)g(z,y)$ on $\Gamma\times\T^n$, which is uniform for $x\in E$. 
\end{proof}

\begin{remark}\label{remark:const-reg-g} It is noteworthy that the effective boundary data inherits the regularity of the given data; i.e., $\overline{g}\in C^k(\Gamma)$ if $g\in C^k(\Gamma\times\T^n)$ for any $k\geq 0$, which is clear from \eqref{eq:gov-const}. We will see later that it fails when interior oscillation comes into the play.
\end{remark}

\begin{remark}\label{remark:const-non-osc} Proposition \ref{proposition:const} can be straightforwardly generalized to non-oscillatory elliptic systems, i.e., $\overline{A} = \overline{A}(x)$ with $x$ being the physical variable on domain $\Omega$.
\end{remark}

\subsection{Boundary Layer Correctors and Periodic Characterization}\label{subsection:bdry-cor}

Let us now consider the case when the interior coefficients of \eqref{eq:Le} periodically oscillate in $\e$-scales. Our first task is to study half-space problems, following the idea of \cite{GM1} and \cite{GM2}.

Let $\chi^*:\T^n\ra \R^{n\times m^2}$ be the matrix corrector for the adjoint $A^*$, i.e., the unique weak solution to 
\begin{equation}\label{eq:chi-pde}
\begin{dcases}
- \nabla \cdot (A^*(y) \nabla (\chi^*(y) + p(y))) =0 & \text{in }\T^n,\\
\int_{\T^n} \chi^*(y) dy = 0,
\end{dcases}
\end{equation}
where $p:\R^n\ra \R^{n\times m^2}$ is the mapping defined by 
\begin{equation}\label{eq:p}
p_i(y) = y_iI\quad (y\in\R^n,1\leq i\leq n).
\end{equation} 
Since $A^*\in C^{1+[\frac{n}{2}]}(\T^n;\R^{n^2\times m^2})$, we deduce from the interior regularity theory for elliptic systems that $\chi^*\in C^{1+ [\frac{n}{2}],\gamma}(\T^n;\R^{n\times m^2})$ with 
\begin{equation}\label{eq:chi-Cka}
\norm{\chi^*}_{C^{1+ [\frac{n}{2}],\gamma}(\T^n;\R^{n\times m^2})} \leq C,
\end{equation}  
where $0<\gamma<1$ is chosen by \eqref{eq:gamma}. Indeed, one may replace $\gamma$ by any positive number less than $1$. 

With the presence of $\chi^*$, let us consider an elliptic system on a half-space, 
\begin{equation}\label{eq:V-pde}
\begin{cases}
-\nabla \cdot (A^* (y) \nabla V^*(y)) = 0 & \text{in } H,\\
V^*(y) = - \chi^*(y) & \text{on } \Pi,
\end{cases}
\end{equation}
where $H := \{y\in\R^n:\nu\cdot y > d\}$ and $\Pi:=\{y\in\R^n:\nu\cdot y =d\}$ with $\nu\in\Ss^{n-1}$ and $d\in\R$. As is shown in \cite{GM2} (see also \cite{P}), there are two ways to express solution $V^*$, one is by the variational formulation and another is by the Poisson integral representation. We shall first take the variational approach, which allows us to describe the quasi-periodicity of $V^*$ on $\Pi$, and then come back to the integral representation, in order to have some refined estimates. 

The idea behind the variational approach is to split the half-space $H$ into a union of hyperplane $\Pi_0+ (t+d)\nu$ for $t > 0$, with $\Pi_0 := \{y\in\R^n: \nu\cdot y =0\}$, and look at the modulo 1 image of $\Pi_0$ in $\T^n$. To be more specific, let us take $\{\nu^1,\cdots,\nu^{n-1}\}\subset\R^n$ as an orthonormal basis for the hyperplane $\Pi_0$, and denote by $N$ the $(n\times(n-1))$-dimensional matrix whose columns are exactly $\nu^1,\cdots,\nu^{n-1}$; apparently, we have $Ny' \in \Pi_0$ for any $y'\in\R^{n-1}$. We remark that the choice of $\{\nu^1,\cdots,\nu^{n-1}\}$ is defined up to a $(n-1)$-dimensional orthonormal transformation. Under these circumstances, let us consider the following system,
\begin{equation}\label{eq:phi-pde}
\begin{cases}
- (NN^T\nabla_\theta + \nu \p_t) \cdot A^*(\theta + t\nu) (NN^T\nabla_\theta + \nu\p_t) \phi^*(\theta,t) = 0 & \text{in }\T^n\times(0,\infty),\\
\phi^*(\theta,0) = -\chi^*(\theta) & \text{on }\T^n.
\end{cases}
\end{equation}

It should be stressed that the system \eqref{eq:phi-pde} is no longer uniformly elliptic, since $NN^T\nabla_\theta + \nu\p_t$ is a $n$-dimensional derivative involving $(n+1)$-dimensional variables. Nonetheless, this system is well-posed, according to \cite[Proposition 2]{GM1}.

\begin{lemma}\label{lemma:phi-C1a} There exists a unique classical solution $\phi^*$ to \eqref{eq:phi-pde}, satisfying $(NN^T\nabla_\theta + \nu\p_t)\phi^* \in W^{1+ [\frac{n}{2}],2}(\T^n\times[0,\infty);\R^{n^2\times m^2})$ with
\begin{equation}\label{eq:phi-Wk2}
\norm{(NN^T\nabla_\theta + \nu\p_t)\phi^*}_{W^{1+ [\frac{n}{2}],2}(\T^n\times[0,\infty);\R^{n^2\times m^2})} \leq C.
\end{equation} 
In particular, $(NN^T\nabla_\theta + \nu\p_t)\phi^* \in C^{1,\gamma}(\T^n\times[0,\infty);\R^{n^2\times m^2})$ with
\begin{equation}\label{eq:phi-C1a}
\norm{(NN^T\nabla_\theta + \nu\p_t) \phi^*}_{C^{1,\gamma}(\T^n\times [0,\infty); \R^{n^2\times m^2})} \leq C.
\end{equation}
\end{lemma}

\begin{proof} The unique existence of $\phi^*$ and the estimate \eqref{eq:phi-Wk2} follow from \cite[Proposition 2]{GM1}, due to \eqref{eq:Cka-A} and \eqref{eq:chi-Cka}. The $C^{1,\gamma}$-regularity of $(NN^T\nabla_\theta + \nu\p_t)\phi^*$ is a direct consequence of the Sobolev embedding $W^{ 1+ [\frac{n}{2}],2} \ra C^{1,\gamma}$, due to the particular choice of $\gamma$ in \eqref{eq:gamma}.
\end{proof}

We obtain a solution to \eqref{eq:V-pde} from Lemma \ref{lemma:phi-C1a}.  

\begin{lemma}\label{lemma:V-C1a} Define $V^*:\overline{H} \ra \R^{n\times m^2}$ by 
\begin{equation}\label{eq:V-phi}
V^*( y ) = \phi^* ( y - (\nu\cdot y)\nu + d\nu, \nu\cdot y - d)\quad\text{on }\overline{H}. 
\end{equation}
Then $V^*$ is a solution to \eqref{eq:V-pde} satisfying $\nabla V^* \in C^{1,\gamma}(\R^n;\R^{n^2\times m^2})$ with 
\begin{equation}\label{eq:V-C1a}
\norm{\nabla V^*}_{C^{1,\gamma}(\R^n;\R^{n^2\times m^2})} \leq C. 
\end{equation}
\end{lemma}

\begin{proof} Invoking the change of variables, $y= (\theta + d\nu) + (t-d)\nu$, with $\theta = y - (\nu\cdot y )\nu$ and $t = \nu\cdot y$, we observe that
\begin{equation}\label{eq:V-phi-nu}
\nabla_y V^* (y) = (NN^T\nabla_\theta + \nu\p_t) \phi^*(y-(\nu\cdot y)\nu + d\nu, \nu\cdot y - d).
\end{equation}
Since $\phi^*$ is a classical $C^2$-solution to \eqref{eq:phi-pde}, we observe that $V^*$ solves \eqref{eq:V-pde}. Moreover, the estimate \eqref{eq:V-C1a} follows immediately from \eqref{eq:phi-C1a}.
\end{proof}

So far we have taken the variational approach to find a solution of \eqref{eq:V-pde}. Although this formulation enables us to describe the oscillatory behavior of $V^*$ in terms of a periodic function, it does not provide the $L^\infty$-estimate; indeed, one may notice that Lemma \ref{lemma:phi-C1a} and Lemma \ref{lemma:V-C1a} only yield uniform estimates for the derivatives of $\phi^*$ and $V^*$. This problem can be amended by considering Poisson integral representation. 

The existence of the Green matrix, hence the Poisson kernel $\overline{P}_H^*$, for elliptic system \eqref{eq:V-pde} in half-space $H$ is proved in \cite{DK} for dimension $n=2$ and in \cite{HK} for dimension $n\geq 3$. The following Poisson integral representation of $V^*$ is due to \cite[Proposition 2.4]{GM2} and \cite[Theorem 3.5]{P}. Notice that such a representation is by no means clear, as an elliptic system on half-space lacks maximum principles, even for scalar case ($m=1$).
 
\begin{lemma}\label{lemma:V-Linf} Let $V^*$ be defined by \eqref{eq:V-phi}. Then 
\begin{equation}\label{eq:V-PH}
V^*(y) = - \int_\Pi P_H^* (y,\xi) \chi^*(\xi) d\xi\quad (y\in H). 
\end{equation}
In particular, there holds 
\begin{equation}\label{eq:V-Linf}
\norm{V^*}_{L^\infty(\overline{H})} \leq C. 
\end{equation}
\end{lemma}

\begin{proof} The integral representation \eqref{eq:V-PH} is proved in \cite[Theorem 3.5]{P} whose proof follows closely that of \cite[Proposition 2.4]{GM2}, provided the solution $V^*$ is $C^\infty$. However, the argument works equally well for classical $C^2$-solutions. As noted in \cite[Proposition 2.4]{GM2}, the $L^\infty$-estimate \eqref{eq:V-Linf} is a simple consequence of \eqref{eq:V-PH}, since $P_H^*$ satisfies
\begin{equation}\label{eq:PH-Linf}
|P_H^*(y,\xi)| \leq \frac{C (\nu\cdot y -d)}{|y - \xi|^n}\quad(y\in H, \xi\in \Pi),
\end{equation}
according to \cite[Lemma 2.5(ii)]{GM2}. 
\end{proof}

\subsection{Oscillatory Elliptic Systems}

As discussed heuristically in deriving \eqref{eq:Pe-peri}, we are going to study the homogenization of Poisson kernels for \eqref{eq:Le}. It is well known (c.f. \cite{BLP}) that the family $\{A(\frac{\cdot}{\e})\}_{\e>0}$ of periodically oscillating coefficients homogenizes to the effective coefficient,
\begin{equation}\label{eq:Aov}
\overline{A} = \int_{\T^n} ( A(y) \nabla ( p(y) + \chi(y))) dy,
\end{equation}
where $\chi$ is the solution to \eqref{eq:chi-pde} with $A^*$ replaced by $A$, and $p$ is defined by \eqref{eq:p}. Note that $\overline{A} \in \R^{n^2\times m^2}$ also satisfies the ellipticity condition \eqref{eq:ellip-A}. In the sequel, let us denote by $\overline{P}$ the Poisson kernel for $-\nabla \cdot (\overline{A}\nabla )$ in $\Omega$. 

For notational convenience, let us define $B:\Gamma\times\T^n\ra \R^{m^2}$ and $\overline{B}:\Gamma\ra \R^{m^2}$ by 
\begin{equation}\label{eq:B-Bov}
B(z,y) = \left[A_{ij}(y)\nu_i(z)\nu_j(z)\right]\quad\text{and respectively }\overline{B}(z) = \left[\overline{A}_{ij} \nu_i(z) \nu_j(z)\right]. 
\end{equation} 
That is, we have $B^{\alpha\beta}(z,y) = A_{ij}^{\alpha\beta}(y)\nu_i(z)\nu_j(z)$ and $\overline{B}^{\alpha\beta}(z)= \overline{A}_{ij}^{\alpha\beta}\nu_i(z)\nu_j(z)$ for $1\leq \alpha,\beta\leq m$. One should note that $\overline{B}(z)$ is an invertible matrix for each $z\in\Gamma$, due to the uniform ellipticity of $\overline{A}$. 

Next we consider the Dirichlet corrector $\Phi_\e^*:\overline\Omega\ra \R^{n\times m^2}$ for $-\nabla\cdot(A^*(\frac{\cdot}{\e})\nabla )$ in $\Omega$, which is given by the solution of 
\begin{equation}\label{eq:Phie-pde}
\begin{dcases}
-\nabla\cdot\left( A^*\left(\frac{x}{\e}\right) \nabla \Phi_\e^*(x)\right) = 0 & \text{in }\Omega,\\
\Phi_\e^*(x) = p(x) & \text{on }\Gamma.
\end{dcases}
\end{equation} 
Let us denote by $P_\e$ the Poisson kernel for $-\nabla\cdot (A(\frac{\cdot}{\e})\nabla)$ in $\Omega$. According to \cite[Theorem 1]{AL2}, the homogenization of $P_\e$ occurs, provided $\Gamma$ is $C^{1,\mu}$; i.e., 
\begin{equation}\label{eq:Pe}
P_\e(x,z) = \overline{P}(x,z)\omega_\e(z) + R_\e(x,z)\quad(x\in\Omega,z\in\Gamma),
\end{equation}
where $\omega_\e:\Gamma \ra \R^{m^2}$ is defined by 
\begin{equation}\label{eq:ome}
\omega_\e (z) = \left[\nu_k(z) \p_\nu \Phi_{\e,k}^*(z)\right] \overline{B}(z)^{-1}  B\left(z,\frac{z}{\e}\right),
\end{equation}
with $\p_\nu$ being the inward normal derivative on $\Gamma$, and the remainder term $R_\e:\Omega\times\Gamma\ra \R^{m^2}$ satisfies 
\begin{equation}\label{eq:Re}
\lim_{\e\ra 0} \sup\{ |R_\e(x,z)|: x\in E,z\in\Gamma\} =0,
\end{equation}
for each $E\Subset\Omega$; although \cite[Theorem 1]{AL2} is stated for the scalar case ($m=1$), the arguments work equally well for systems. Let us also remark that a quantitative estimate on the decay rate of $R_\e$ is established in \cite[Theorem 1.1]{KLS} for dimension $n\geq 3$, provided $\Gamma$ is $C^{2,\mu}$. However, we will only use the qualitative result \cite[Theorem 1]{AL2} in the sequel, which holds for any dimension $n\geq 2$ and any $C^{1,\mu}$-regular boundary $\Gamma$. 

In what follows, we describe the oscillatory behavior of $\omega_\e$ in terms of a periodic mapping, based on the idea of \cite{AKMP}. Due to $C^{1+[\frac{n}{2}]}$-regularity of $A$, we are able to define $\phi^*:\Gamma\times \T^n\times [0,\infty) \ra \R^{n\times m^2}$ such that
\begin{equation*}
\text{$\phi^*(z,\cdot,\cdot)$ is the solution of \eqref{eq:phi-pde} with $\nu=\nu(z)$}\quad(z\in\Gamma),
\end{equation*}
according to Lemma \ref{lemma:phi-C1a}. We observe from \eqref{eq:phi-C1a} and the fact $N(z)^T\nu(z) = 0$ that $\p_t\phi^*(z,\cdot,\cdot)\in C^{1,\gamma}(\T^n\times[0,\infty);\R^{n\times m^2})$ uniformly for $z\in\Gamma$ with
\begin{equation}\label{eq:phiz-C1a}
\norm{\p_t\phi^*(z,\cdot,\cdot)}_{C^{1,\gamma}(\T^n\times[0,\infty);\R^{n\times m^2})} \leq C.
\end{equation} 

With $\phi^*$ at hand, we define $w:\Gamma\times\T^n\ra \R^{m^2}$ by 
\begin{equation}\label{eq:w}
w (z,y) := \left[I + \nu_k(z)(\p_t \phi_k^*(z,y,0) + \p_\nu \chi_k^*(y))\right] \overline{B}(z)^{-1} B(z,y),
\end{equation}
where $\p_\nu$ is again the inward normal derivative on $\Gamma$. Then it is clear from \eqref{eq:ellip-A}, \eqref{eq:phiz-C1a} and \eqref{eq:chi-Cka} that $w\in L^\infty(\Gamma\times\T^n;\R^{m^2})$, and that $w(z,\cdot)\in C^{1,\gamma}(\T^n;\R^{m^2})$ uniformly for $z\in\Gamma$ with
\begin{equation}\label{eq:w-C1a}
\norm{w(z,\cdot)}_{C^{1,\gamma}(\T^n;\R^{m^2})} \leq C.
\end{equation}

The next lemma essentially follows the proof of \cite[Corollary 5.2]{AKMP}. The estimate \eqref{eq:ome-w} becomes less sharp, however, due to the lack of convexity of the domain $\Omega$. Let us provide the proof for the sake of completeness. 

\begin{lemma}\label{lemma:ome-w} For $\omega_\e$ and $w$ defined by \eqref{eq:ome} and \eqref{eq:w} respectively, there are universal constants $C$ and $r_0$ such that 
\begin{equation}\label{eq:ome-w}
\left|\omega_\e(z) - w\left(z,\frac{z}{\e} \right) \right|\leq C\left( \frac{\e}{r} + r^\gamma + \frac{r^{(1+\mu)}}{\e} \right),
\end{equation}
for any $z\in\Gamma$, any $\e>0$ and any $0< r \leq r_0$.
\end{lemma}

\begin{proof} Let us fix $z\in\Gamma$, $\e>0$ and $0<r\leq r_0$, where $r_0$ will be chosen later. In comparison of $\omega_\e$ with $w$, we observe that the crucial step of this proof is to show
\begin{equation}\label{eq:Phie-C01}
\left|\p_\nu \Phi_\e^* (z) - \p_\nu p(z) - \p_t \phi^*\left(z,\frac{z}{\e},0\right) - \p_\nu \chi^*\left(\frac{x}{\e}\right)\right| \leq C\left(\frac{\e}{r}+r^\gamma + \frac{r^{(1+\mu)}}{\e} \right),
\end{equation}
for certain universal $C>0$. The idea in showing \eqref{eq:Phie-C01} is to approximate the system \eqref{eq:Phie-pde} in a local neighborhood $B_{2r}(z)\cap\overline\Omega$ by a half-space problem of type \eqref{eq:V-pde}. 

Let us choose a half-space which contains the local neighborhood so that the associated half-space solution is well-defined on $B_{2r}(z)\cap\overline\Omega$. Recall from \eqref{eq:H} the definition of the half-space $H(z)$, whose boundary is precisely the tangent hyperplane $\Pi(z)$ to $\Gamma$. Owing to the $C^{1,\mu}$-regularity \eqref{eq:C1a-bdry} of $\Gamma$, we can choose a universal $r_0$ satisfying 
\begin{equation}\label{eq:C1a-bdry-re}
\overline\Omega\cap B_{2r}(z)\subset H(z) - Kr^{1+\mu}\nu(z)\quad\text{for $0<r\leq r_0$}.
\end{equation}

Let us define $V^*$ by \eqref{eq:V-phi} with $\phi^* = \phi^*(z,\cdot,\cdot)$, $\nu=\nu(z)$, $d = \e^{-1}(\nu(z)\cdot z - Kr^{1+\mu})$ and $H=\e^{-1} ( H(z) - Kr^{1+\mu}\nu(z))$. Then by the construction of $\phi^*$ and Lemma \ref{lemma:V-C1a}, $V^*$ is a solution to \eqref{eq:V-pde}. Hence, it follows from \eqref{eq:Phie-pde} and \eqref{eq:chi-pde} that the mapping $v_\e: \overline{\Omega}\cap B_{2r}(z) \ra \R^{n\times m^2}$, defined by  
\begin{equation}\label{eq:ve}
v_\e(x) := \Phi_\e^*(x) - p(x) - \e V^*\left(\frac{x}{\e}\right) - \e\chi^* \left(\frac{x}{\e}\right),
\end{equation}
solves the elliptic system
\begin{equation}\label{eq:ve-pde}
\begin{dcases}
- \nabla \left( A^* \left(\frac{x}{\e}\right) \nabla v_\e(x) \right) = 0 & \text{in } \Omega\cap B_{2r}(z), \\
v_\e(x) = - \e V^*\left(\frac{x}{\e}\right) - \e\chi^* \left(\frac{x}{\e}\right) & \text{on } \Gamma\cap B_{2r}(z).
\end{dcases}
\end{equation}

Due to the H\"{o}lder regularity of the coefficient $A^*$ followed by \eqref{eq:Cka-A}, the uniform boundary Lipschitz estimate \cite[Lemma 20]{AL1} applies to \eqref{eq:ve}, which yields that 
\begin{equation}\label{eq:ve-C01}
\begin{split}
&\norm{\nabla v_\e}_{L^\infty(\Omega\cap B_r(z);\R^{n^2\times m^2})} \\
&\leq C \left\{ \frac{1}{r} \norm{v_\e}_{L^\infty(\Omega\cap B_{2r}(z);\R^{n\times m^2})} + r^\gamma \left[ \nabla \left( ( \e V^* + \e \chi^*)\left(\frac{\cdot}{\e}\right) \right)\right]_{C^\gamma(\Gamma\cap B_{2r}(z);\R^{n^2\times m^2})} \right\}.
\end{split}
\end{equation}
The $L^\infty$-estimate for $v_\e$ can be derived as follows. Since $V^* + \chi^* = 0$ on $\Pi = \e^{-1}(\Pi(z) - Kr^{1+\mu}\nu(z))$, we may continuously extend $V^*+\chi^*$ by zero in $\R^n\setminus\overline{H}$; here the continuity across $\Pi$ is ensured by from \eqref{eq:chi-Cka} and \eqref{eq:V-C1a}. Then $V^* + \chi^*$ is well-defined on $\e^{-1}\Gamma$, and by \eqref{eq:chi-Cka} and \eqref{eq:V-Linf} we have $\norm{V^* + \chi^*}_{L^\infty(\e^{-1}\Gamma;\R^{n\times m^2})} \leq C$. Incorporating the uniform $L^\infty$-estimate \cite[Theorem 3(i)]{AL1} for the Poisson kernel $P_\e^*$ for $-\nabla \cdot (A^*(\frac{\cdot}{\e}) \nabla)$ in $\Omega$, i.e., $|P_\e^*(x,z)| \leq C\dist(x,\Gamma)|x-z|^{-n}$ for any $x\in\Omega$ and any $z\in\Gamma$, we derive that
\begin{equation}\label{eq:ve-Linf}
\begin{split}
|v_\e(x)| & \leq \left| \int_\Gamma P_\e^* (x,z)\left(\e V^*\left(\frac{z}{\e}\right) + \e \chi^*\left(\frac{z}{\e}\right)\right) d\sigma_z\right| \leq C\e \norm{V^*+ \chi^*}_{L^\infty(\e^{-1}\Gamma)} \leq C\e,
\end{split}
\end{equation}
for all $x\in\Omega$. On the other hand, we deduce from \eqref{eq:chi-Cka} and \eqref{eq:V-C1a} that 
\begin{equation}\label{eq:ve-C1a}
\left[ \nabla \left( ( \e V^*+ \e \chi^*)\left(\frac{\cdot}{\e}\right) \right)\right]_{C^\gamma(\Gamma\cap B_{2r}(z);\R^{n^2\times m^2})} \leq C,
\end{equation}
since $\Gamma\cap B_{2r}(z) \subset H$ by \eqref{eq:C1a-bdry-re}. Combining \eqref{eq:ve-Linf} and \eqref{eq:ve-C1a}, we arrive at 
\begin{equation}\label{eq:ve-C01-re}
\norm{\nabla v_\e}_{L^\infty(\Omega\cap B_r(z);\R^{n^2\times m^2})}\leq C\left(\frac{\e}{r} + r^\gamma\right).
\end{equation}

We observe from \eqref{eq:V-phi-nu} that
\begin{equation}\label{eq:Vezr-phiz-nu}
\p_\nu V^* \left(\frac{z}{\e}\right) = \p_t \phi^*\left( z, \frac{z - Kr^{1+\mu}\nu(z)}{\e}, \frac{Kr^{1+\mu}}{\e}\right),
\end{equation}
where $\p_\nu$ is the directional derivative in $\nu(z)$. By means of the Lipschitz estimate of $\p_t\phi^*(z,\cdot,\cdot)$ deduced from \eqref{eq:phiz-C1a}, we have 
\begin{equation*}
\left| \p_t \phi^*\left(z,\frac{z - Kr^{1+\mu}\nu(z)}{\e}, \frac{Kr^{1+\mu}}{\e}\right) - \p_t \phi^* \left(z,\frac{z}{\e},0\right) \right| \leq \frac{CKr^{1+\mu}}{\e}.
\end{equation*}
This inequality combined with \eqref{eq:ve-C01-re} yields \eqref{eq:Phie-C01}.
\end{proof}

As a corollary, we have the following homogenization of the Poisson kernel $P_\e$. 

\begin{corollary}\label{corollary:Poisson} For $w$ defined by \eqref{eq:w}, there holds 
\begin{equation}\label{eq:Pe-Pov-w}
P_\e(x,z) = \overline{P}(x,z) w\left(z,\frac{z}{\e}\right) + R_\e(x,z) \quad(x\in\Omega,z\in\Gamma),
\end{equation}
where $R_\e:\Omega\times\T^n\ra \R^{m^2}$ is the remainder term satisfying \eqref{eq:Re} for each $E\Subset\Omega$. 
\end{corollary}

\begin{proof} Let us choose $r=\e^\kappa$ for any $(1+\mu)^{-1} <\kappa <1$ in \eqref{eq:ome-w} and insert it to \eqref{eq:Pe}. The rest of the proof then follows easily from the fact that the Poisson kernel $\overline{P}$ satisfies $|\overline{P}(x,z)| \leq C\dist(x,\Gamma)|x-z|^{-n}$ for any $x\in\Omega$ and $z\in\Gamma$, as shown in \cite[Theorem 3(i)]{AL1}.
\end{proof}

We are now ready to prove the boundary layer homogenization of \eqref{eq:Le}.

\begin{proof}[Proof of Theorem \ref{theorem:bdry-hom}]
Fix $E\Subset\Omega$. Due to Proposition \ref{corollary:Poisson}, we have, uniformly for $x\in E$,  
\begin{equation*}
\lim_{\e\ra 0} u_\e(x) = \lim_{\e\ra 0}\int_\Gamma \overline{P}(x,z) w\left(z,\frac{z}{\e}\right) g\left(z,\frac{z}{\e}\right) d\sigma_z.
\end{equation*}
Utilizing \eqref{eq:equi-g} and \eqref{eq:w-C1a}, the rest of the proof follows the same argument as in the proof of Proposition \ref{proposition:const}, whence we omit the details. 
\end{proof}

\subsection{Sharpness of IDDC and Counterexamples}

Here we provide some counterexamples for Theorem \ref{theorem:osc-int} and Theorem \ref{theorem:bdry-hom} when IDDC fails for $\Gamma$, proving the sharpness of our results. 

For simplicity, let us consider scalar case with constant coefficients, i.e., \eqref{eq:Love} with $m=1$. The first example treats the case when $\Omega$ is given by a strip normal to the standard basis vector $e^n$, which is apparently a rational direction. We will see that due to the absence of IDDC, we are able to generate all the possible limits within the maximum and minimum values of $g$. 

\begin{example}\label{example:limit} Let $\overline{A}$ be an identity matrix in $\R^{n^2}$, $\Omega = \{x\in\R^n: |x_n|<R\}$ with a certain $R>0$, and take a function $h\in C^\infty(\T)$ with a constant $M$ such that $g(R,y) = h(y_n)$ on $\T^n$ and that $g(-R,y) = M$ on $\T^n$. Notice that $\Gamma$ is a union of two hyperplanes $x_n=R$ and $x_n=-R$, so it does not satisfy IDDC. 

Note that the modulo $1$ image of $\{\e^{-1}R\}_{\e>0}$ cover the entire space $\T$. This implies that there exists a sequence $\e_i\ra 0$ such that $h(\e_i^{-1}R) = d$ for all $i=1,2,\cdots$. In this way, the solution $\overline{u}_{\e_i}$ to \eqref{eq:Love} satisfies $\overline{u}_{\e_i}(x) = M$ for $x_n=-R$ and $\overline{u}_{\e_i}(x) = d$ for $x_n=R$, for all $i$. As $-\nabla \cdot(\overline{A} \nabla)$ being Laplacian, $\overline{u}_{\e_i}(x) = \frac{d-M}{2R}(x_n+R) + M$ in $\Omega$. Hence, the limit profile $\overline{u}$ becomes the same affine function, which is precisely the solution to \eqref{eq:Lov} with $\overline{g}(x) = M$ if $x_n=-R$ and $\overline{g}(x) = d$ if $x_n=R$. This shows that limit profiles of $\{\overline{u}_\e\}_{\e>0}$ may vary depending on the choice of the sequence $\e=\e_i$, and attain all possible values in $h(\T)$ on the hyperplane $x_n=R$. 
\end{example}

Next we see that the effective boundary data $\overline{g}$ may turn out to be discontinuous, if $\Gamma$ contains a flat piece normal to a rational direction.

\begin{example}\label{example:discont} Let $\overline{A}$ be the identity matrix in $\R^{n^2}$, $\Omega = \{x\in\R^n: |x_n|<R, |x'|< 1 + (R^2-x_n^2)^{1/2}\}$ for some $R>0$, and $g$ a smooth periodic function such that $g(x,y) = h(y_n)$ for some non-constant $h\in C^\infty(\T)$. Note that $\Gamma$ is the union of $\Gamma_1 := \{x\in\R^n: |x_n|<R, |x'| = 1+ (R^2 - x_n^2)^{1/2}\}$ and $\Gamma_2 := \{x\in\R^n: |x_n|=R\}$. $\Gamma_1$ is a union of two hemispheres, whence satisfies IDDC, while $\Gamma_2$ is a flat piece with a rational normal. 

Let us write by $\overline{h}$ the average value of $h$, i.e., $\overline{h} = \int_\T h(t)dt$. Note that $\overline{h} < \max_\T h$,  since $h$ is not constant. Then any possible limit solution $\overline{u}$ of $\{\overline{u}_\e\}_{\e>0}$ assumes value $\overline{h}$ on $\Gamma_1$. However, as we take a sequence $\e_i\ra 0$ such that $h(\e_i^{-1}R) = \max_\T h$, we have $\overline{u}(x) = \lim_{i\ra \infty} \overline{u}_{\e_i}(x) = \max_\T h$ on $\Gamma_2$. Thus, the boundary data is not continuous on $\overline\Gamma_1\cap\Gamma_2$. 
\end{example}

\end{document}